\documentclass{article} 
\begin{document}
 
  \begin {center} 

{\bf  SEMIGROUP IDENTITIES, PROOFS, AND ARTIFICIAL INTELLIGENCE  }\\

Sherman Stein

Mathematics Department

Mathematics Sciences Building

University of California at Davis

1 Shields Avenue

Davis, CA 95616

e-mail: stein@math.ucdavis.edu

key words: identities, models, proofs, computer proofs, Burnside's Problem

  \end {center}

\newpage

\title {SEMIGROUP IDENTITIES AND ARTIFICIAL INTELLIGENCE  }
\author{ Sherman Stein}


\begin{center}  
{SEMIGROUP IDENTITIES AND ARTIFICIAL INTELLIGENCE }\\
 Sherman Stein 

\end{center}

(Abstract)  It is known that if every group     satisfying an identity of the form $yx \sim xU(x,y)y$ is abelian, 
so is every semigroup
that satisfies that identity.  
Because a group has an identity element and the cancellation property, it is easier to show that a group is abelian than that a semigroup is. If we know that it is, then there must be a sequence of substitutions using $xU(x,y)y \sim yx$ that transforms $xy$ to $yx.$ We examine such sequences and propose finding them as a challenge to proof by computer.
\vspace{.1cm}

Also, every model of $y \sim xU(x,y)x$ is a group. This raises a similar challenge, which we explore in the special case $y \sim x^my^px^n.$ In addition we determine the free model with two generators of some of these identities. In particular, we find that the free model for $y \sim x^2yx^2$ has order $32$ and is the product of $D_4$ (the symmetries of a square), $C_2$, and $C_2$, and point out relations between such identities and Burnside's Problem concerning models of $x^n= e$.






\vskip.3in

\section { Introduction.}

 This paper concerns identities in two variables on
semigroups, their models, and their implications. They are of interest for three
reasons:
\vspace{.1cm}

$\bullet $They raise the question, Is there an algorithm for deciding whether a given
identity is a consequence of another identity?
\vspace{.1cm}
 
$\bullet $They offer a source of challenges to artificial intelligence to find shorter proofs
of the theorems here or proofs of other identities. (One identity is the axiom and the other is the alleged  theorem.)
\vspace{.1cm}

$\bullet $They provide a different perspective on some groups.

\section {Terminology.}

    A $semigroup$ is a set with an associative binary
operation. A $word$ in $x$ and $y$ is a finite string of $x$s and $y$s. Its
length is the number of letters in it. For instance, the word $xyxyyy$ has
length six. We will use exponential notation, writing such a word as $xyxy^3. $
  An $identity$ in $x$ and $y$ consists of two words set equivalent, $W_1(x, y) \sim   W_2(x, y).$    A semigroup satisfies that identity if whenever $x$ and $y$ are replaced by elements of the semigroup a true statement results. For instance, a semigroup that satisfies the identity
$xy \sim   yx$ is commutative. A semigroup that satisfies the identity $x\sim y$ has
only one element.    A semigroup that satisfies a given identity is a $model$ for
that identity. 
\vspace{.2cm}
 
 If all semigroups that satisfy an identity $I$ also satisfy an identity $J$, then $J$ can be deduced from $I$ by a sequence of substitutions using $I$.
If one reverses the order of the symbols in $I$ and $J$ to produce the identities $I'$ and $J'$, then $J'$ is a consequence of $I'$. We call $I'$ and $J'$ the duals of $I$ and $J$ and $J' \sim I'$ the dual identity of the identity $I \sim J.$
\vspace{.2cm}

Let $U=U(x,y)$ be a word in $x$ and $y$ and $P$ and $Q$ also words in $x$ and $y$. Replacing $x$ by $P$ and $y$ by $Q$ throughout $U$ produces a word which we denote $U(P,Q).$ Let $W_1(x,y) \sim W_2(x,y)$ be the axiom. We will say two words are equivalent if one can be written as $RW_1(P,Q)S$ and the other as $RW_2(P,Q)S$ for some words $P,Q,R,S$ in $x$ and $y$. This generates an equivalence relation which we denote $\sim .$  If $U(x,y) \sim V(x,y)$ it follows that $U(P,Q) \sim V(P,Q)$ for any two words $P$ and $Q$ in $x$ and $y$. To show this, consider the sequence of substitutions using the axiom that transforms $U(x,y)$ into $V(x,y)$. Throughout the words that appear in this sequence replace $x$ by $P$ and $y$ by $Q$. This produces a sequence of substitutions using the axiom that tranforms $U(P,Q)$ into $V(P,Q)$.
\vspace{.1cm}

We label a proof using group theory G and a proof using only substitutions S for "semigroup" or "substitution." An S-proof modeled on a G-proof we will label GS. If two words in an S-proof are identical we use the equal sign =. If they are obtained by using one or more substitutions, we write $\sim.$
\vspace{.1cm}
\vspace{.2cm}

\section{Identities $yx \sim xUy$ and commutativity.}

An identity $yx\sim xUy$ has the form 
$$yx \sim x^{a_1}y^{b_1}x^{a_2}y^{b_2} \dots x^{a_n} y^{b_n}, \eqno(1)$$
where the $2n$ exponents are positive. We will examine several instances of (1).
\vspace{.1cm}

The following theorem due to Tully treats the case $n=1$ of (1) [4].
\vspace{.1cm}

{\bf Theorem 3.1.} For positive integers $p$ and $q,$ $yx \sim x^p y^q$ implies $yx \sim xy.$
 \vspace{.1cm}
 
 {\it G-Proof.} Replacing $x$ by the identity element $e$ shows that $y=y^q$. Similarly, $x = x^p$. Thus $x^p y^q = xy,$ which implies $yx= xy$.\hfill {\tt[]}
\vspace{.1cm}

In an S-proof the parentheses in a word enclose the active section, where the substitution is applied. 
\vspace{.1cm}

 {\it  S-Proof.} The case when $pq = 1$ is trivial, so assume that $pq \geq 2.$ Because
$x^{p+q}\sim  x^2$, exponents can be altered by $t=p+q-2.$    Thus, if $r\geq 2$, then $x^r\sim  x^{r+t}$. The following sequence establishes 
commutativity:
$$yx\sim  x^{p}y^{q}\sim  (y^{pq})(x^{pq})\sim  (y^{t}y^{pq})(x^{pq}x^{t})\sim $$
$$ y^{p^2q}(y^{qt}x^{pt}x^{pq^2 })\sim  (y^{p^2q}x^{pq^2})\sim  x^{pq}y^{pq}\sim  y^{p}x^{q}\sim  xy.$$
    \hfill {\tt[]}

The next theorem concerns an identity that appeared as a problem [2].
\vspace{.1cm}

{\bf Theorem 3.2.} For a positive integer $n$, $yx \sim (xy)^n $ implies $yx \sim xy.$
\vspace{.1cm}

{\it G-Proof.} Replacing $y$ by the identity element shows that $x\sim x^n.$ Thus $(xy)^n \sim xy$. Hence $yx \sim xy.$
                                                                                              {\tt[]}
\vspace{.1cm}

{\it S-Proof.}

$$ yx\sim (xy)^n \sim (yx)^{n^2} = [(yx)(yx)^{n-1}]^n \sim(yx)^{n-1}yx = (yx)^n \sim xy.$$\hfill {\tt[]}

We turn next to identities of the form $$yx \sim x^{a_1}y^{b_1}x^{a_2}y^{b_2}, \eqno(2)$$
the case $n=2$ of (1).
\vspace{.1cm}

When all four exponents in (2) are $1$ we have a special case of Theorem 3.2. So we turn to the case when exactly one exponent is greater than $1$, namely $yx \sim x^nyxy$ and $yx \sim xy^nxy$, the other two cases being duals of these.
\vspace{.1cm}

 In the proofs we will use the identities to change exponents. For instance,
replacing $y$ by $x$ shows that $x^2 \sim x^{a_{1}+b_{1}+a_{2}+b_{2}}$. Thus exponents can be changed by $a_{1}+b_{1}+a_{2}+b_{2}-2$ as long as they remain at least $2$. Replacing $y$ by $x^2$ and then $x$ by $y^2$ shows that exponents can also be changed by $a_{1}+2b_{1}+a_{2}+2b_{2}-3$ and by 
$2a_{1}+b_{1}+2a_{2}+b_{2}-3$ as long as they remain at least $3$. This implies that exponents can be changed by $a_1 +a_2 -1$ and by $b_1 +b_2 -1$.
\vspace{.1cm}

{\bf Theorem 3.3.} If exactly one of the exponents is greater than $1$, then the identity (2) implies $yx \sim xy.$
\vspace{.1cm}

{\it G-Proof.} Consider $yx \sim x y^n x y$ for $n\geq 2.$ Replacing $y$ by $e$ shows that $x \sim x^2$ and cancellation shows that $e = x.$ Thus any group satisfying (2) has only one element and is therefore commutative. A similar argument goes through for $yx\sim x^nyxy.$ {\tt[]}
\vspace{.1cm}

{\it S-Proof.} The identity $yx \sim xy^nxy$ implies $yx \sim xy^2xy$ by exponent changes. Specifically, replacing $y$ by $x$ in the identity shows that $x^2\sim x^{n+3}$, implying that exponents can be raised or lowered by $n+1$ as long as they remain at least $2$.  Replacing $x$ by $y^2$ shows that $y^3\sim y^{n+5}$, implying that exponents can be changed by $n+2$ as long as they remain at least $3$. If $n\geq 3$, we may raise $n$ by $3((n+5)-3)$ to reach $4n+6$, and then reduce $4n+6$ by $4((n+3)-2)$ to reach $2.$ Thus it suffces to prove the theorem in the special case $n=2, yx\sim xy^2xy.$ We then have
$$yx \sim xy^2(xy) \sim x(y^2y)x^2yx \sim xy^2x^2yx = xy(yx^2yx) \sim xyxy.$$
Thus $yx \sim xyxy$, a case  considered in Theorem 3.2.{\tt[]}
\vspace{.1cm}

Next we consider the case when exactly two of the exponents in (2) are at least $2$. There are four types and their duals:
$$(i)  yx\sim x^ay^bxy , (ii) 
 yx\sim x^ayx^by, (iii) yx \sim xy^ax^by, (iv)  yx \sim x^ayx y^b.$$
\vspace{.1cm}

{\bf Theorem 3.4.} Let $a$ and $b$ be at least $2$. Then (i), (ii), and (iii) imply commutativity, but (iv) does not.

{\it G-Proof.} For (i) replace $x$ and then $y$ by $e$ showing that $y^{b+1}\sim y$ and $x^{a+1} \sim x$. Cancellation then gives $y^b=e$ and $x^a=e$,  Thus (i) implies $yx \sim xy.$ Case (iii) is similar.
\vspace{.1cm}

For (ii) set $x$ equal to $e$, obtaining $y^2 = y$. Cancellation gives $y=e$, and the group has only one element, so is commutative.
\vspace{.1cm}

For (iv), if $(a,b) \geq 3$, then any group for which $x^{(a,b)} = e$ satisfies the identity. Therefore there are nonabelian groups that satisfy the identity. If $(a,b) = 1$ or $2$, then any group model satisfies the identity $x^2 = e$, hence is abelian.{\tt[]}

\vspace{.1cm}

{\it S-Proof.} Case (i). Let $d = (a,b).$ For $d\geq 2$ exponents can be altered by $d$ as long as they stay at least $d$. Thus we can assume
$yx \sim x^d y^d xy.$ Then we have
$$yx \sim x^dy^d(xy) \sim (x^dy^d)y^dx^dyx \sim (x^{d^2}y^{d^2}y^dx^d)yx \sim(y^dx^d)yx \sim xy.$$
If $d=1$, we can begin with $yx \sim x^2y^2xy$ instead.
\vspace{.1cm}

Case (ii) Exponents can be altered by $a+b$ or $2(a+b)-1$ as long as they stay at least $2$. Because $(a+b,2(a+b)-1)=1$ they can be altered by $1$ as long as they stay at least $2$. Therefore it suffices to consider only the special case $yx \sim x^2yx^2y.$: 
$$yx \sim (x^2y)x^2y \sim y^2x^2y^2(x^2x^2)y \sim(y^2)x^2(y^2)x^2y$$ 
$$\sim (y^4x^2y^4x^2)y \sim x^2y^2y = x^2(y^3) \sim x^2y^2.$$
Thus
$$yx \sim x^2y^2 \sim y^4x^4 \sim y^2x^2.$$
Finally,
$$yx \sim x^2y^2 \sim y^2x^2 \sim xy.$$

Case (iii) Exponents can be altered by $a+b,a+2b,$ and $2a+b,$ hence by $a$, as long as they stay at least $2$. In particular $b$ can be changed to $a$ because $b + 2(2a+b)-3(a+b)=a.$ Therefore we need show only that $yx\sim xy^ax^ay$ implies $yx\sim xy.$ We have
$$yx\sim (xy^a)x^ay\sim y^ax^a(y^{a^2}xx^ay\sim y y^ax^ay^axx^ay=y^a(x^ay^a)x^axy\sim (y^ax^{a^2}y^{a^2}x^a)xy\sim x^ay^axy.$$
Thus $yx\sim x^ay^axy$, which is case (i). {\tt[]}
\vspace{.1cm}

 Case (iv) Any group that satisfies $x^{(a,b)}=e$ satisfies the identity. Since there are nonabelian groups that satisfy $x^{(a,b)}=e$ for $(a.b)\geq 3$, the identity does not imply commutativity. However, if $a,b)$ is $1$ or $2$, it does imply commutativity. We take the case $(a,b) = 1$ and leave the case $(a,b)=2$ to the reader. 
\vspace{.1cm}

Exponents can be changed by $1$. We have
$$yx \sim (x^2y)xy^2 \sim (y^2x^2)yx^4xy^2 $$
 $$\sim x^4y^2x^2(y^4yx^4x)y^2 \sim x^4(y^2)x^2(y^4)x^4y^2 $$

$$\sim x^4(y^4x^2y^2x^4)y^2 \sim (x^4x^2)(y^2y^2) \sim x^2y^2.$$
Thus
$$yx\sim x^2y^2 \sim y^4x^4 \sim y^2x^2\sim xy.$$                                                           
\hfill  {\tt[]}

The next five theorems sample some of the possibilities when three or all of the four exponents are at least $2.$
\vspace{.1cm}

 {\bf Theorem 3.5.} If $a,b,$ and $c$ are at least $2$, the identity $yx \sim  xy^ax^by^c$ implies $yx \sim xy$.
\vspace{.2cm}

{\it G-Proof.} Replacing $x$ and $y$ by $e$ shows that $x^b= e$ and $y^{a+c} = y$. Thus $yx=xy^{a+c} = xy.$ {\tt[]}
\vspace{.1cm}

The reader may provide an S-proof.
\vspace{.1cm}

On the other hand, the identity $yx\sim x^ayx^by^c$ may or may not imply $yx\sim xy,$ as the next two theorems show.

{\bf Theorem 3.6.} If $a,b,$ and $c$ are at least $2$, $d=(a,c) \geq 3$ and $b\equiv 1$ (mod $d$) then the identity $yx\sim x^ayx^by^c$ does not imply 
$yx\sim xy$.

{\it G-Proof.} A nonabelian group satisfying $x^d=e$ also has $x^b = x,x^a = e=y^c,$ and therefore satisfies $yx=xy$. {\tt[]}
\vspace{.1cm}

In the opposite direction we have
\vspace{.1cm}

{\bf Theorem 3.7.}If $a,b,$ and $c$ are at least $2$ and $(a+b-1,c)$ is $1$ or $2$, the identity $yx\sim x^ayx^by^c$ implies $yx\sim xy.$
\vspace{.1cm}

{\it G-Proof.} Replacing $x$ and $y$ by $e$ shows that $y^c=e$ and $x^{a+b}=x$, hence $x^{a+b-1}=e.$ Thus $x^{(c,a+b-1)}=e$, which implies that $yx\sim xy.$ {\tt[]}
\vspace{.1cm}

When all four exponents are at least two we again address only a few cases.
\vspace{.1cm}

  {\bf Theorem 3.8.} If $a,b,c$ and $d$ are at least $2$ and $(a+c-1,b+d-1)$ is $1$ or $2$, then the identity $yx\sim x^ay^bx^cy^d$ implies $yx\sim xy$.

\vspace{.1cm}

 {\it G-Proof.} Set $x$ and $y$ equal to $e.$ {\tt[]}

However, even if $(a+c-1,b+d-1)$ is greater than $2$, the identity can still imply $yx\sim xy$, as the following example and Theorem 3.9 show.
 Consider, for instance, $yx \sim x^5y^6x^4y^3$, for which
$(a+c-1,b+d-1)=8$. In a group, because $x=x^9,$ exponents can be changed by $8,$ and we have 
$$ yx \sim  x^5 (y^6 x^4) y^3 \sim x^5(x^4)^5 (y^6)^6 (x^4)^4 ( y^6)^3 y^3 \sim  x^{25} y^{36} x^{16} y^{21} \sim x y^{4} e y^{5}  \sim x y^9 \sim xy.$$

{\bf Theorem 3.9.}  Let $k$ be at least $2$. The identity $yx\sim x^ky^kx^ky^k$ implies $yx\sim xy$ if and only if $k \equiv 0$ or $1$ (mod $3$.)
\vspace{.1cm}

{\it G-Proof.} Consider first $k\equiv 2$ (mod $3$). A group for which $x^3=e$ satisfies the identity if and only if it satisfies the identity $yx\sim  x^2y^2x^2y^2.$ As noted in [5], it does satisfy the latter identity for  $yx \sim x^2y^2x^2y^2\sim x^{-1}y^{-1}x^{-1}y^{-1}$, which is equivalent
 to $(yx)^3=e.$
\vspace{.1cm}

The case $k \equiv 0$ or $1$ depends on the fact that $x^{6}y^{6}\sim y^{6}x^{6}$, which holds for any value of $k$. Note that exponents can be changed by $4k-2.$ We have
$$x^6y^6=x^2(x^4y^6)\sim x^{4k}(x^4y^6)\sim (x^{4k}y^{6k}x^{4k}y^{6k})x^{4k}\sim (y^6x^4 x^{4k})\sim y^6x^4x^2 =y^6x^6.$$

\vspace{.1cm}

To show that if $k\equiv 0$ or $1$ (mod $3$ ) then $yx\sim xy$, we first alter exponents to be a multiple of $6$ to exploit the weak form of commutativity just obtained.
Therefore, we wish to show that there is a positive integer $n$ such that $k\equiv 6n$ (mod $2k-1$).(Exponents can be changed by $4k-2$ and $6k-3$, hence by $2k-1.$) This amounts to showing that there are integers $n$ and $u$ such that $(2k-1)u=k-6n.$ We may assume that if there is a solution, then there is a solution with $n$ positive.

Because $k\equiv 0$ or $1$ (mod $3$) we have $(6,2k-1)=1$. Hence the equation $k \equiv 6n$ (mod $2k-1$) is solvable. Thus
$$yx\sim x^{6n}y^{6n}x^{6n}y^{6n}=(x^6)^n(y^6)^n(x^6)^n(y^6)^n.$$
Because $x^6$ and $y^6$ commute,
$$yx\sim (y^6)^n(x^6)^n(y^6)^n(x^6)^n\sim xy.
\hfill $$ {\tt[]}

\section {\bf Identities $y\sim xUx$.}
 Let $U$ be a word in $x$ and $y$. For any model of $y \sim xUx$ the equations $ax=b$ and $ya=b$ can be solved for $x$ and $y$. Thus any model is a group, hence has an identity element, each element has an inverse, and  left- and right-cancellations hold. There are therefore two ways to establish that $y\sim   xUx$ implies an identity $I$: either use properties of a group or use only substitutions based on the identity $y \sim xUx$. Some examples will illustrate the two approaches.
\vspace{.2cm}

{\bf Theorem 4.1.} For a positive integer $n$, $y \sim (xy)^nx$ implies $x^2\sim y^2$ and $xy\sim yx.$  If $n$ is even, $y\sim (xy)^nx$ implies $x\sim y.$
\vspace{.1cm}

{\it G-Proof.} Replace $x$ by $e$, obtaining $y^n\sim y.$ Thus $xyx\sim y$. Replacing $y$ by $e$ in $xyx\sim y$ shows that $x^2=e$, hence $x^2\sim y^2$ and $xy\sim yx$.

For any $n$ we have $x^{n+1}=e$ and $y^{n-1}=e$, thus $x^{(n+1,n-1)}=e$. If $n$ is even $(n+1,n-1)=1$ and we conclude that $x=e$, thus $x\sim y.$ {\tt[]}
\vspace{.1cm}

We give two S-proofs, the first using only one exponent change, the second using two.
\vspace{.1cm}

{\it S-Proof.} We have $$y^2= y(y)\sim y(xy)^nx=((yx)^ny)x \sim x^2.$$
Then we use the exponential change $y \sim y^{2n+1}$. Thus
$$xy \sim (xy)^{2n} xy = x(yx)^n ((yx)^n y)$$
$$\sim x(yx)^nx=((xy)^nx)x \sim yx.$$
If $n$ is even, we have
$$y \sim [(xy)^2]^{n/2} x \sim [(x^2)^2]^{n/2}x = (xx)^n x \sim x.$$ \hfill  {\tt []}

{\it S-Proof.} Replacing $x$ by $y$ and by $y^2$ shows that exponents can by changed by $2n$ and $3n+1$. Thus the exponent $n$ can be changed to $1$, for $n +(3n+1) - 2(2n)=1.$ Consequently $y \sim xyx.$ Hence
$$y^2 = y(y) \sim (yxy)x \sim xx = x^2.$$ Also we have $(x)(y) \sim (yx\cdot yx\cdot yx)\sim yx.$

 If $n$ is even $(2n,3n+1)=1$, so $n$ can be reduced to $2$. We then have
$$y \sim ((xy)^2)x \sim x^2 x=xxx\sim x.$$ \hfill {\tt[]}

Some proofs consist of a sequence of substitutions that do not decrease the lengths of the intermediate words, followed by a sequence that does not increase the lengths. They amount to two reductions of some longer word. Such a proof is called a "proof of degree one" or a "mountain proof" in [4].  It is shown there that for the identity $xyx\sim y$ any proof can be replaced by a mountain proof but that is not the case for $(xy)^2\sim y^3x^3.$
We offer a mountain proof showing that $y\sim (xy)^2x$ implies $x\sim y$, an instance of Theorem 4.1.
\vspace{.1cm}

 With $x$ replaced by $xyx^2$ and $y$ replaced by $x$,we have
$$x\sim(xyx^2)x(xyx^2)x(xyx^2)=xyx(x^2)xyx(x^2)xyx^2$$
$$\sim xyxyx^2yx^2yxyx(yx^2yx^2y)xyx^2$$
$$=xy(xy\cdot x\cdot xy\cdot x\cdot xy)xyx(yx\cdot x\cdot yx\cdot x\cdot yx)yx^2$$
$$\sim xy(x\cdot xy\cdot x\cdot xy\cdot x)x$$
$$\sim xy\cdot xy\cdot x\sim y$$ \hfill {\tt[]}

As Martin Davis has suggested, one may obtain an S-proof by imitating a G-proof. We
 illustrate this approach by proving the case $n = 2$ of Theorem 4.1, namely $y\sim (xy)^2x$ implies $x\sim y$, the one just obtained by a mountain S-proof.
\vspace{.1cm}

 {\it GS-Proof.} Because $x\sim  x^5$ we expect $x^4$ to behave like the identity element.  Indeed it does, for $x^4x\sim x$  and $x^4y\sim x^4(xyxyx)=x^5yxyx\sim xyxyx\sim y.$ Similarly $yx^4\sim y.$ Also $y^4$ is an identity element. Thus $x^4\sim x^4y^4=y^4$, and so $x^4\sim y^4.$ We have

$$y\sim x^4yx^4yx^4\sim x^{12}y^2=(x^4)^3y^2\sim y^2.$$
From $y\sim y^2$ it follows that $y^2\sim y^4$, hence $y\sim y^4\sim x^4\sim x.$ {\tt[]}
\vspace{.1cm}

 Now we turn our attention to identities of the form $x^my^px^n=y$, where $m,n,$ and $p$ are positive integers. 
 \vspace{.2cm}

{\bf Theorem 4.2.} The identity $ y\sim x^my^px^n $ implies $y\sim x^myx^n,y^p\sim y,x^{m+p+n}\sim x,x^{m+p+n-1}y\sim y\sim yx^{m+p+n-1},
x^my\sim   yx^m,x^ny\sim yx^n,x^{m+n}\sim y^{m+n}$.
\vspace{.1cm}

 {\it G-Proof.} We have $y\sim (x^2)^my^p(x^2)^n=x^m(x^m y^p x^n)x^n\sim x^m y x^n.$ Replacing $x$ by $e$ in the given identity shows that $y^p = y.$ Replacing $y$ by $x$ instead shows that $x^{m+p+n} =x.$ Next, $$x^{m+n+p-1}y=x^{m+n+p-1}x^my^px^n =x^{m+n+p}x^{m-1}y^px^n=xx^{m-1}y^px^n=x^my^px^n\sim y.$$
 Similarly, $yx^{m+n+p-1}=y.$ From $x^m x^n=e$ we have $x^n = x^{-m}.$ Thus $x^myx^n=y$ implies $x^myx^{-m}=y$ or $x^my=yx^m$. Similarly $x^ny=yx^n.$ Thus $x^{m+n}y^p x^{m+n}=
x^n(x^my^px^n)x^m=x^nyx^m=x^myx^n=y.$ Finally, from $x^{m+n}=e$ follows $x^{m+n}=y^{m+n}.$         {\tt[]}
\vspace{.1cm}

{\it S-Proof.} The identity $x^myx^n\sim y$ is proved the same way as in the G-proof. Exponents can be altered by $p-1$ and $m+n$ as long as they remain at least one. Thus $y^p \sim y.$ The identities $x^{m+n+p}\sim x$ and $x^{m+n+p-1}y\sim yx^{m+n+p-1}$ are proved as in the G-proof. From $x^myx^n\sim y$ follows
$x^myx^n x^{m+p-1}\sim yx^{m+p-1}.$ Thus $x^my \sim yx^{m+p-1}= yx^{m-1+p} \sim yx^{m-1+1} = yx^{m}.$ Similarly $x^ny \sim yx^n$ and therefore 
$x^{m+n}y^px^{m+n}=y.$ Finally, $x^{m+n}=(x^m)(x^n)\sim y^m(x^my^ny^mx^n)y^n \sim x^my^ny^mx^n \sim y^{m+n}.$    {\tt[]}
\vspace{.2cm}
 
{\bf Theorem 4.3.} The identity $y\sim x^my^px^n$ implies $y^{p-1}\sim x^{p-1}$
\vspace{.1cm}

{\it G-Proof.} By Theorem 4.2, $y^p=y$. Cancellation gives $y^{p-1}=e$, which implies $y^{p-1}\sim x^{p-1}.$  {\tt[]}

\vspace{.2cm}

We offer three S-proofs. The first, by Martin Davis, a GS-proof, makes use of an identity element. The second, an S-proof, consists mostly of a proof of degree one. The third, due to Dean Hickerson, is a proof of degree one. 
\vspace{.1cm}

{\it S-Proof.} From $x^{m+n+p} \sim x$ we have $x^{m+n+p-1}x\sim x$, which suggests that $x^{m+n+p-1}$ is an identity element. We show that it is by obtaining $yx^{m+n+p-1} \sim y:$
$$(y)x^{m+n+p-1}\sim x^my^px^nx^{m+n+p-1}=x^my^px^{n-1}( x^{m+n+p}) \sim x^myx^{n-1}x = x^myx^n \sim y.$$

 A similar argument shows that  $x^{m+n+p-1}y \sim y$.

\vspace{.1cm}

Let $Q=x^{m+n+p-1}$. Then $y^p \sim Q^m y^p Q^n \sim y.$ Thus
$$y^{p-1}(x) \sim y^{p-1}y^mx^py^n = (y^p) y^{m-1} x^py^n \sim y y^{m-1} x^py^n = y^mx^py^n \sim x.$$

 Similarly $yx^{p-1} \sim y.$ Thus $y^{p-1} \sim y^{p-1}x^{p-1} \sim x^{p-1}.$ {\tt[]}

\vspace{.1cm}

 {\it S-Proof.} The preceding remarks show that there is an integer $d$ such that $y\sim x^dy^px^d$. Working with this identity, we have
$$y^{p-1}\sim(y^dx^p)^d y^{(p-1)p} (y^dx^p)^d.\eqno (3)$$

 If $d$ is even, we have
$$y^{p-1}\sim [(y^dx^py^d)x^p]^{d/2} y^{(p-1)p} [(y^dx^py^d)x^p]^{d/2}.$$

 By an exponent change, $y^{(p-1)p}y^d$ reduces to $y^d.$ Thus
$$y^{p-1} \sim(x\cdot x^p)^{d/2}(x\cdot x^p)^{d/2}= x^{(p+1)d}.$$
 Replacing $y$ by $x$ shows that $x^{p-1}\sim x^{(p+1)d}$, hence $x^{p-1} \sim y^{p-1}$. 
\vspace{.1cm}

 If $d$ is odd, (3) is
$$y^{p-1} \sim [(y^dx^py^d)x^p]^{(d-1)/2} y^dx^px^{(p(p-1)}y^dx^p [(y^d x^py^d)x^p]^{(d-1)/2}.$$

An exponent change removes $y^{p(p-1)}$ and we have
$$y^{p-1} \sim (x\cdot x^p)^{(d-1)/2}x \cdot x^p (x \cdot x^p)x^{(d-1)/2} = (x \cdot x^p)^d.$$
Thus $y^{p-1} \sim x^{p-1}.$   {\tt[]}
\vspace{.2cm}

The following S-proof of the same theorem collects several shorter proofs into one mountain proof.
\vspace{.1cm}

{\bf S-Proof.}  Let $Z$ denote $y^{p+1}x^{p+1}.$ Then
$$y^{p-1} \sim Z^m y^{p(p-1)}Z^n =Z^{m-1} y^{p+1}(x^{p+1})y^{p^2 + 1} x^{p+1}Z^{n-1}$$  
$$\sim Z^{m-1} y^{p+1}x^m x^{p(p+1)} x^n y^{p^2 +1} x^{p+1} Z^{n-1}$$  
$$=Z^{m-1}y^{p+1} x^{p^2+m+p+n} y^{p^2}(y) x^{p+1} Z^{n-1}$$  
$$\sim Z^{m-1} y^{p+1} x^{p^2 +m+p+n} y^{p^2}y^m y^p y^n x^{p+1}Z^{n-1}$$  
$$=Z^{m-1}y^{p+1}x^{p^2 + m+p+n} y^{p^2 + m +p +n} x^ {p+1} Z^{n-1}$$  
$$=Z^{m-1} y^{p+1}x^{p^2}(x^mx^px^n) y^{p^2 +m+p+n} x^{p+1} Z^{n-1}$$ 

$$\sim Z^{m-1}y^{p+1}x^{p^2} x y^{p^2 +m+p+n} x^{p+1} Z^{n-1}$$  
$$=Z^{m-1}y^{p+1}x^{p^2 +1}(y^m y^{p(p+1)} y^n) x^{p+1}Z^{n-1}$$  
$$ \sim Z^{m-1}y^{p+1} x^{p^2+1} y^{p+1} x^{p+1} Z^{n-1}$$ 

$$=Z^m  x^{p(p-1)}Z^n \sim x^{p-1}.$$ \hfill {\tt []}

\section {\bf Other Identities.} 

So far we have considered only identities related to groups. We now examine some identities not related to them.
\vspace{.1cm}

For a word $Q$ let $f(Q)$ be the number of ordered pairs, called "forward pairs," $x$ and $y$ with $x$ to the left of $y$. For instance, $f(x^3y^3xy)=13.$ Let $g(Q)$ be the number of $x$s in $Q$ and $h(Q)$ be the number of $y$s in $Q$. For words $U,V,W,W'$ such that $g(W)=g(W')$ and $h(W)=h(W')$ one has $f(UWV)-f(UW' V)= f(W)-f(W').$
This equation will be used several times in this section.

\vspace{.1cm}

The first part of the proof of the next theorem is due to Dean Hickerson.
\vspace{.1cm}

{\bf Theorem 5.1.} Let $a$ and $b$ be positive integers, $a<b.$ Then $(xy)^a\sim (yx)^b$ implies $(xy)^a\sim (yx)^n$ if and only if $b-a$ divides $n-a.$

\vspace{.1cm}

{\it S-Proof.} We record the word $(xy)^p(yx)^q$ by the ordered pair $(p,q)$. Then we have
$$(0,a)\sim (b,0)\sim (b-a,b)\sim (2b-a,b-a)\sim (2b-2a,2b-a) \sim ...$$

Thus $(0,a)\sim (k(b-a),a+k(b-a))$ for $k\geq 0.$
\vspace{.1cm}

Similarly,
$$(a,0)\sim (0,b)\sim (b,b-a)\sim (b-a,2b-a)\sim (2b-a,2b-2a)\sim (2b-2a,3b-2a)\sim ...$$

Thus $(a,0) \sim (k(b-a),b+k(b-a))$ for $k\geq 0.$
\vspace{.1cm}

If $a+k(b-a)\geq 2a$ an exponent change can increase it to $b+k(b-a)$. Specifically, from $x^{2a}\sim  x^{2b}$ and $x^{3a}\sim  x^{3b}$ exponents can be changed by $b-a=2(2b-2a)-(3b-3a).$ Thus $(k(b-a),a+k(b-a))\sim (k(b-a),b+k(b-a),$ from which it follows that $(0,a)\sim (a,0).$
At this point we have $(xy)^a\sim (yx)^a.$ Thus $(yx)^a\sim (yx)^b.$ If $n-a$ is divisible by $b-a$ $n$ can be written as $a+k(b-a)$ for some non-negative integer $k$. Therefore
$$(yx)^n =(yx)^{a+kb-ka}\sim (yx)^{a+ka-ka} = (yx)^a \sim (xy)^a.$$
Conversely, if $(xy)^a\sim (yx)^b$, then $b-a$ divides $n-a$. This follows from the fact that $g(UV)^a)-g(UV)^b$ is a multiple of $b-a$.

 \hfill {\tt[]}

The preceding proof that $(xy)^a\sim (yx)^a$ is a mountain proof. The inequality $a+k(b-a)\geq a$ is equivalent to $k\geq a/(b-a).$ If $a/(b-a)$ is large, then the number of steps would be large. The longest word, the peak of the proof, has length $2k(b-a)+b$, which is at least $2a+b$ and at most $3b.$ There may be proofs using fewer steps or shorter words.

\vspace{.2cm}

{\bf Theorem 5.2} Let $a$ and $b$ be positive integers, $a<b.$ Then $(xy)^a \sim (yx)^b$ implies $(xy)^a \sim (xy)^n$ if and only if $b-a$ divides $n-a$. 

\vspace{.1cm}

{\it S-proof.} By Theorem 5.1, $(xy)^a\sim (yx)^a$. But $(yx)^a\sim (xy)^b$; hence $(xy)^a\sim (xy)^b.$ It follws that if $b-a$ divides $n-a$ then 
$(xy)^a \sim (xy)^n/$ If $(xy)^a\sim (xy)^n$, then $b-a$ divides $n-a$ by the same reasoning as in the proof of Theorem 5.1.
\hfill  {\tt []}
 
\vspace{.1cm}

{\bf Theorem 5.3.} Let $m$ and $n$ be positive integers. Then $(xy)^m \sim (yx)^m$ implies $(xy)^n \sim (yx)^n$ if and only if $m$ divides $n$.
\vspace{.1cm}

{\it S-proof.} Assume that $m$ divides $n$, $n = km.$ Then
$$(xy)^n=[(xy)^m]^k \sim [(yx)^m]^k = (yx)^n.$$
\\ To treat the case when $m$ does not divide $n$, we use the function $f$ that records the number of forward pairs. We have $f((xy)^n) = \Sigma_{i=1}^ni$ and
$f((yx)^n) = \Sigma_{i= 1}^{n-1}i.$ Their difference is $n$. As may be checked, for any words $U$ and $V$, $f((UV)^m)- f((VU)^m) $ is a multiple of 
$m$. Thus $m$ must divide $n$.             {\tt[]}
\\

{\bf Theorem 5.4.} Let $m$ and $n$ be positive integers.  Then $xy^m \sim y^m x$ implies $(xy)^n \sim (yx)^n$ if and only if $m$ divides $n$, and $n>m$. 
\vspace{.1cm}
\\ 

{\it S-Proof.} On each use of the identity $xy^m \sim y^mx$ the number of forward pairs changes by a multiple of $m$, as may be checked. However,
$f((xy)^n)-f((yx)^n) = n.$ Thus $n$ must be a multiple of $m$.
 Now assume that $n>m$ and there is an integer $k$ such that $n = km.$ Then $k \geq 2$ and

$$(xy)^n = (xy)^{mk}=[(xy)^m]^k=x[(yx)^m]^{k-1}(yx)^{m-1}y$$
$$ \sim[(yx)^m]^{k-1}x (yx)^{m-1}y$$
$$=[(yx)^m]^{k-2}(yx)^{m-1}(yx)x(yx)^{m-1}y$$
$$=[(yx)^m]^{k-2}(yx)^{m-1}y(x(xy)^m)$$
$$\sim [(yx)^m]^{k-2}(yx)^{m-1}(yx)^{m+1}$$
$$=(yx)^{km} = (yx)^n.$$
The case $n=m$ is not a consequence of the identity $xy^m\sim y^mx$ because $(xy)^m$ has no section of the form $UV^m$.   
\hfill {\tt[]}
\vspace{.2cm}

A similar argument establishes the next theorem.
\vspace{.1cm}

{\bf Theorem 5.5.} If $x^my^n\sim y^nx^m$ implies $x^ay^b\sim y^bx^a$, then $mn$ divides $ab.$
\vspace{.1cm}

The converse does not hold. For instance, $x^my^n\sim y^nx^m$ does not imply $x^ay^{mn}\sim y^{mn}x^a$ if $a$ is not a sum of $m$s and $n$s. On the other hand, we
have the following theorem.
\vspace{.1cm}

{\bf Theorem 5.6.} The identity $x^my^m\sim y^mx^m$ implies $x^ay^a\sim y^ax^a$ if and only if $m$ divides $a$.
\vspace{.1cm}

{\it Proof.}  By Theorem 5.5, $m^2$ divides $a^2$, hence $m$ divides $a$. Also, if $m$ divides $a$, $x^ay^a\sim y^ax^a$. {\tt[]}
\vspace{.1cm}

A direct proof of Theorem 5.6 follows from the fact that if $a$ is not a multiple of $m$, then any word equivalent to $x^ay^a$ begins with $x$.
\vspace{.1cm}

{\bf Theorem 5.7} The identity $xy^m\sim y^mx$ implies $xy^n\sim y^nx$ if and only if $m$ divides $n$.
\vspace{.1cm}

The proof of this is similar to the preceding proofs.

\section{\bf Free Models}
Let $I$ be an identity in $x$ and $y$ and $a$ and $b$ be two symbols. The words in $a$ and $b$ form a semigroup with the multiplication being juxtaposition. Set two words equivalent if the identity obtained by replacing $a$ and $b$ by $x$ and $y$ respectively is a consequence of $I$. The equivalence classes form the free semigroup on two letters for $I$, which we denote $F(I)$. For instance, $F(x\sim y)$ has only one element and 
$F(xy \sim yx)$ can be viewed as the abelian semigroup consisting of the elements $a^ib^j$, where $i,j \geq 0$ and $i+j >0,$  with the product of $a^ib^j$ and $^{i'}b^{j'}$ being $a^{i+i'}b^{j+j'}.$ We shall examine $F(x^my^px^n \sim y)$ for certain values of $m,n,$ and $p$.  
\vspace{.1cm}

By Theorem 4.2, if $(m,n) = 2$,  $x^2$ any square can be moved past any part of a word. This is the basis of the following lemma.
The symbol $a^0$ or $b^0$ indicate the absence of $a$ or $b$ respectively.
\vspace{.1cm}

{\bf Lemma. 6.1}    Let $I$ be an identity of the form $x^m y x^n \sim   y$ with $(m,n) = 2.$ Then every element of $F(I)$ has a representative of the form
$a^i b^j , a^i b^j a$, or $a^i b^j ab,$ where $i, j, \geq 0$.
\vspace{.1cm}

 {\it Proof.}    Consider a word $a^{i_1} b^{j_1} a^{i_2} b^{j_2} \ \cdots \ a^{i_m} b^{j_m}.$ Because any square can be moved to the left, we can assume that the
word is equivalent to one of the form $a^i b^j abab. . . ab(a), i,j\geq 0$, two
exponentials followed by a tail of alternating $a$s and $b$s.
\vspace{.1cm}

Next we shrink the tail to one of at most two letters. To do this, assume that the length of the tail is at least three and begins
 $a^i b^j aba.$ If $j$ is
$0$, this is $a^{(i+1)}ba.$ If $j$ is positive and even, $b^j$ can be moved past
$a$ which can merge with $a^i$, resulting in $a^{(i+1)} b^{(j+1)}a$. If $j$ is odd and $i$ is positive, we have this sequence of
transformations:
                $$a^i b^j aba = a^i b^{(j-1)} baba \sim   b^{(j-1)} a^i baba =$$
$$ b^{(j-1)} a^{(i-1)}ababa \sim  b^{(j-1)} a^{(i-1)} a^2 bab = b^{(j-1)} a^{(i+1)} bab \sim   a^{(i+1)} b^j ab.$$   

The case when $j$ is odd and $i$ is $0$ can be treated in fewer steps. Thus every element of the free semigroup 
for the identity has a representative in which
the tail has at most two elements.     {\tt  []}
\vspace{.1cm}

Each element of $F(x^m y x^n \sim   y)$ satisfies the equation $x^{(m+n)} = e$.   
This implies that the exponents $i$ and $j$ in the lemma can be chosen less than
$m + n$. (For book-keeping purposes we may use $m + n$ as a substitute for the
exponent $0.$) Consequently the free model is finite, having at most $3(m +
n)^2$ elements. Because the number of $a$s and also the number of $b$s changes
by a multiple of    $m + n$ on each use of the axiom, the elements represented by $a^i b^j, 0 \leq i, j < m + n$ are not equal. Therefore the order of 
$F(x^myx^n \sim y)$ is at least $(m + n)^2$.    In case $(m,n) = 2$ and $m + n$ is a multiple of $4$, the order is $2(m+n)^2$, as the next theorem shows.
\vspace{.1cm}   
           
The proof of the theorem depends on Lemma 6.2, which concerns the
number of forward pairs.    
   \\ 
         
 {\bf Lemma 6.2.} Assume that $m$ and $n$ are even and that $m + n$ is a
multiple of $4.$    Then on each substitution using the identity $x^m y x^n \sim y$ the parity of the number of forward pairs does not change.

{\it Proof.}    Let $U, V, W,$ and $T$ be words in $a$ and $b$. Assume
that the word $UVW$ is transformed by the substitution of $T^m V T^n$ for $V,$
resulting in the word $UT^mVT^nW$. Assume that $U$ has $u\  a$s and $u' \ b$s, $V$ has $v$ $a$s and $v' \  b$s, $W$ has $w\  a$s and $w' \  b$s, and $T$ has $t\  a$s and $t'\  b$s. Also assume that $T$ has $f$ forward pairs.    The difference in the number of forward pairs between the initial word $UVW$ and the word obtained by the substitution is

$$(m+n)f +mut' +nu t' +(m+n)(m+n-1)t t'/2+ mt(v'+ w')+ n t w' +n v t'.$$

That sum is even. If the substitution is used in the reverse direction, the
change in the number of forwards pairs is again even, being equal to the above
sum.

This establishes the lemma.{\tt[]}
\vspace{.1cm}

We use the two lemmas to determine the distinct elements of $F(x^2 y x^2 \sim y).$

First of all, each element of that group is represented by a word of the
form

$$ a^i b^j \ {\rm or}\ a^i b^j a  \ {\rm or}\  a^i b^j ab, {\rm with}\ 0 \leq i, j \leq 3.  $$

There are three words of these forms that have $i\ a$s and $j\ b$s,
namely,

$$ a^i b^j , a^{(i-1)}b^j a ,\ {\rm and}\ a^{(i-1)} b^{(j-1)} ab. \eqno(4)$$

The number of forward pairs in each, is respectively,

$$ ij, ij - j, \ {\rm and}\    ij-  j + 1. $$
   
We show that exactly two of these three words are equivalent.
 \vspace{.1cm} 
          
Take the case when $j$ is even. Then $ij$ and $ij-j$ are both even and $ij-j+1$ is odd.    Note that in this case $a^{(i-1)} b^j a$  can be transformed into $a^i b^j$ by moving $b^j$ to the left and then back all the way to the right.
The case when $j$ is odd breaks into two cases depending on whether $i$ is odd
or even. As may be checked, the two words with the same parity are equivalent.
\vspace{.1cm}

For a given number of $a$s and a given number of $b$s, there are
therefore exactly two inequivalent words of the form (1). Thus $F(x^2 y x^2 \sim  
y)$ has $2(4^2) = 32$ elements. Because $ab$ is not equivalent to $ba,$ the
group is not abelian.    We list its elements by the number of $a$s and the number
of $b$s modulo $4,$ with the first word having an even number of forward pairs,
the second, an odd number. The couplet $ij$ accompanies the two words with $i
\ a$s and $j\ b$s. The exponents $0$ and $4$ are interchangeable.

$$ \pmatrix {
  00 & a^4 & a^3b^3ab & 20 & a^2 & ab^3ab \cr
  01 & b & a^3ba & 21 & a^2b & aba \cr
 02 & b^2 & a^3bab & 22 & a^2b^2 & abab \cr
  03 & b^3 & a^3b^3a & 23 & a^2b^3 & ab^3a \cr
  10 & a & b^3ab & 30 & a^3 & a^2b^3ab    \cr
  11 & ba & ab & 31 & a^2ba & a^3b \cr
  12 & ab^2 & bab & 32 & a^3b^2    &a^2bab \cr
  13 & b^3a & ab^3 & 33 & a^2b^3a & a^3b^3 \cr
}
$$

To determine whether two words $U$ and $V$ are equivalent one need only count
the number of $a$s modulo $4$ the number of $b$s modulo $4$, and find the parity of the number of forward pairs in
each. They are equivalent if and only if they have the same triplet of numbers.
To obtain a semigroup transformation of one into the other using the identity
$x^2yx^2 \sim   y$ one could first transform each into a word on the list and then
piece the two sequences obtained to form a transformation of one word into the
other. However, there well may be shorter proofs.
   \\          

Note that using the description of the elements in terms of the number of $a$s,
the number of $b$s, and the parity of the number of forward pairs, one could
describe several subgroups. For instance the elements for which the total number
of $a$s and    $b$s is even and the number of forward pairs is even forms a
subgroup of order $8$.
   \\
          
The set of squares consists of $a^2, b^2, abab,$ and $a^4$, four elements which
do not form a group.    The center has those elements and $a^2b^2,ab^3ab,a^3bab,$
and $a^3b^3ab.$    It consists of the elements in which both $a$ and $b$ appear an
even number of times. It is also the set of elements of order $1$ or $2$.
Because there are eight of them there are $24$ elements of order $4$. This
information pinpoints the group designated 32/8 in [6], which is the product
$D_4 \times  C_2   \times  C_2,$ where $D_4$ is the dihedral group of order $8$ and $C_2$ is
the cyclic group of order $2$.
    \\            

The same reasoning shows that $F(x^6 y x^ {10} \sim   y)$ has order $2(16)^2 = 512$ and is nonabelian.
$ F = F(x^3 y x^3 \sim   y)$ is quite different. The identity amounts to the two
conditions: $ x^6 = e$ and the cube of every element lies in the center.    The tail in this case consists of a sequence of $a$s, $b$s, $a^2$s and $b^2$s. Because there are endless sequences using only three symbols without a square, there are
arbitrarily long tails without cubes. So the approach in the preceding proof
does not go through here. However $F$ is finite because a group with two
generators such that $x^6 = e$ is known to be finite [1].
   \\          

In that vein, consider $F(x^a y^{(a+1)} x^a \sim   y)$ for any positive integer
$a$.    Such a group is simply a group for which $x^a = e$.    Burnside in 1902
asked whether such a group if finitely generated is necessarily finite. 
The answer is yes for $a = 1,2,3,4,$ and $6$, unknown for $a=5$ and a finite number of values of $a$
and no for an infinity of values of $a$ [1].
   \\

\section{Questions.}

 1. The free semigroup on three generators $a, b,$ and $c$ for the identity $x^2 y x^2 \sim   y$ is finite because any group model satisfying
 $ x^4 = e$ is finite. What is the order of that
semigroup? The analog of the tail that was used in Sec.6 is a string made of
$a, b, $ and $c$ instead of just $a$ and $b.$ Such strings can be arbitrarily long without any squares, that is, without a section of the form $ww,$ where $w$ is a word. Using the tools in Sec. 6, one can show only that its order is at least $3(4^3) = 192.$
   \\          

2. Which identities $I$ in this paper have the property that if $I$ implies identity $J$, then there is a mountain proof of $J$ based on $I$?
\vspace{.1cm}

3. Exponent change plays a big role in the proofs in this paper. Can it always be avoided?
\vspace{.1cm}

4. Complete the analysis of identities of the form $yx\sim x^{a_1}y^{b_1}x^{a_2}y^{b_2}.$ What is the situation for $yx\sim x^{a_1}y^{b_1}x^{a_2}y^{b_2}x^{a_3}y^{b_3}?$
\vspace{.1cm}

5. Is there an algorithm for deciding whether one identity in $x$ and $y$ implies another such identity? If there is none, then some proofs must contain words extremely long when compared to the words in the given identities.
\vspace{.1cm}

6. Is there an algorithm for deciding whether an identity in $x$ and $y$ implies $yx\sim xy$? (It is easy to show that
such an identity must have the form $yx\sim xUy$ or $y \sim xUx$.)
\vspace{.1cm}         

7. The identity $x^4y^9x^4\sim   y$ is equivalent to the assertion that $x^8 = e$
and $x^4$ is in the center of $F(x^4y^9x^4\sim  y)$. Is this group infinite or
finite? If it is infinite, that would imply that a group generated by two elements and satisfying $x^8=e$ is infinite. This suggests
using identities of the form $x^ay^bx^c \sim   y$ or other identities to investigate
Burnside's  problem. 
\vspace{.1cm}

8. Which groups generated by two elements are the free model of some identity in $x$ and $y$?
\vspace{.1cm}

9. How should the complexity of a proof be measured? By the number of substitutions used? By the length of the longest word that appears in it? By the number of peaks in it? "Peak" is defined in [4].
\vspace{.1cm}

10. Is there a shorter proof of Theorem 5.1? 

\vspace{2cm}

{\bf References}
           
\vspace{1cm}

1. S. I. Adian, The Burnside Problem and Identities in Groups, Ergebnisse der
Mathematik und ihrer Grenzgebiete 95, New York, 1979.
\vspace{.1cm}

2. D.MacHale and Michael Oseacoid, Problem E 3109, American Mathematical Monthly 92 (1985), 591
      and solution ibid 94 (1987) 468-469.
\vspace{.1cm}

3. S. K. Stein, Semigroups Satisfying $xy = yg(x,y)x,$ Semigroup Forum 36
(1987), 249-251.

 \vspace{.1cm}

4. S. Stein, The Combinatorial Degrees of Proofs and Equations, Algebra
Universalis 35 (1996), 472-484.

\vspace{.1cm}

5 T. Tamura, Semigroups Satisfying Identity $xy = f(x,y),$ Pac. J. Math 31
(1969), 523-521.

 \vspace{.1cm}

6. A. D. Thomas and G.V. Wood, Group Tables, Orpington, UK, 1980.

\vfill

\end{document}